\documentclass{article}

\usepackage{arxiv}

\usepackage[utf8]{inputenc} 
\usepackage[T1]{fontenc}    
\usepackage{hyperref}       
\usepackage{url}            
\usepackage{booktabs}       
\usepackage{amsfonts}       
\usepackage{nicefrac}       
\usepackage{microtype}      
\usepackage{lipsum}		
\usepackage{graphicx}
\usepackage{natbib}
\usepackage{doi}
\usepackage{graphicx}
\usepackage[linesnumbered,ruled,vlined]{algorithm2e}
\usepackage{amssymb}
\usepackage{lscape} 
\usepackage{appendix}
\usepackage{lineno}
\usepackage{amsfonts}
\usepackage{booktabs}
\usepackage{siunitx}
\usepackage[normalem]{ulem}
\usepackage{hyperref}
\hypersetup{colorlinks=true,linkcolor=blue,filecolor=magenta,  urlcolor=cyan}
\usepackage{caption}
\usepackage{subcaption}
\usepackage[utf8]{inputenc}
\usepackage[T1]{fontenc}
\usepackage{textcomp}
\usepackage[dvipsnames]{xcolor}
\usepackage{amsmath}
\usepackage{amsthm} 
\usepackage{latexsym, amssymb}
\usepackage{array}
\usepackage{xspace}
\usepackage{amsmath, xparse}
\usepackage{multirow}

\title{Enriched finite element approach for modeling discontinuous electric field in multimaterial problems}

\author{{Christian Narváez-Muñoz} \\
        Escola Tècnica Superior d'Enginyers de Camins, Canals i Ports\\
        Universitat Politècnica de Catalunya - BarcelonaTech (UPC)\\
	Centre Internacional de Mètodes Numérics en Enginyeria  (CIMNE)\\
        Barcelona, Spain\\
	\texttt{cnarvaez@cimne.upc.edu} \\
	\And
	{Mohammad R. Hashemi} \\
	Centre Internacional de Mètodes Numérics en Enginyeria  (CIMNE)\\
        Barcelona, Spain\\
	\AND
	{Pavel B. Ryzhakov} \\
        Escola Tècnica Superior d'Enginyers de Camins, Canals i Ports\\
        Universitat Politècnica de Catalunya - BarcelonaTech (UPC)\\
	Centre Internacional de Mètodes Numérics en Enginyeria  (CIMNE)\\
        Barcelona, Spain\\
	\texttt{pavel.ryzhakov@upc.edu} \\
        \AND
	{Jordi Pons-Prats} \\
        Department of Physics, Aeronautics Division\\
        Universitat Politècnica de Catalunya, Barcelona Tech (UPC)\\
	Centre Internacional de Mètodes Numérics en Enginyeria  (CIMNE)\\
        Barcelona, Spain\\
        \AND
	{Herbert Owen} \\
	Department of Computer Applications in Science and Engineering\\
        Barcelona Supercomputing Center (BSC-CNS)\\
        Barcelona, Spain
}

\renewcommand{\shorttitle}{\textit{arXiv} Template}

\hypersetup{
pdftitle={A template for the arxiv style},
pdfsubject={q-bio.NC, q-bio.QM},
pdfauthor={David S.~Hippocampus, Elias D.~Striatum},
pdfkeywords={First keyword, Second keyword, More},
}

\begin{document}
\maketitle

\begin{abstract}
This work is devoted to the development of an efficient and robust technique for accurate capturing of electric field in multi-material problems. 
The formulation is based on the finite element method enriched by the introduction of hat-type shape function within the elements crossed by the material interface.
The peculiar feature of the proposed method consists in direct employment of the hat-function that requires solely one additional degree of freedom per cut element for capturing the discontinuity in electric potential gradient and, thus, the electric field.
This additional degree of freedom is subsequently statically condensed element-wise prior to the assembly of the global discrete system. As a consequence, the graph of the system matrix remains the same as that of the standard finite element method.
In order to guarantee the robust performance of the proposed method for a wide range of electrical material property ratios, it also accounts for the possible discontinuities among the neighboring cut elements that arises due to employing additional degrees of freedom fully local to the element. The method is tested using several examples solved on structured and unstructured grids. The proposed approach constitutes a basis for enriched FEM applicable to a wide range of electro-magnetic problems.
\end{abstract}


\section{Introduction}
In electromagnetic problems involving multi-materials, one faces field discontinuities that pose challenges for the numerical modeling. For instance, in electro-mechanical (EME)~\cite{rochus2011electrostatic} and electro-hydrodynamics (EHD)~\cite{narvaez2021determination} systems, discontinuity in the permittivity (in case of dielectrics) or conductivity(in case of conductors) at the material interface results in a discontinuous gradient of the electric potential, and, consequently the electric field therein~\cite{ramos2011electrokinetics}. Being based on a weak formulation, which naturally incorporates the gradient boundary conditions, the finite element method (FEM) has been by far the most commonly used approach to treat the problems in the area of electro-statics/dynamics~\cite{hansbo2004finite,moes2017extended}. In case of multi-materials, FEM allows accurate representation of the discontinuity at the materials' interface in case of using interface-conforming meshes~\cite{hansbo2004finite}. However, ensuring mesh conformance is not straight-forward in the problems where the material interface constantly evolves and/or deforms. This occurs, for instance, in the field of EHD~\cite{segura2020sketch,tomar2007two}, where liquid undergoes continuous deformation and/or motion under the effect of the electric field. For such problems, in case of using fixed grid finite elements in standard form, as soon as the material interface moves away from the mesh nodes, a significant error in the gradient of the potential is introduced as the discontinuous gradient becomes represented by a continuous one with a sharp slope. Such errors immediately compromise the quality of the numerical solution since the net electric force acting at the material interface is directly proportional to the jump in the gradient of the electric potential~\cite{melcher1969electrohydrodynamics}. Thereby, it is essential to enhance FEM in order to improve the accuracy of the overall model with respect to discontinuity capturing.\\

Different numerical FEM-based strategies have been devised in the literature for capturing the discontinuities without the need to use a conforming mesh. The generalized or extended Finite Element Method (henceforth in this study we refer to it as G/X-FEM), provides a valuable contribution in handling such problems~\cite{melenk1996partition, belytschko1999elastic, moes1999finite}. This is achieved by enriching the solution space via introducing additional degrees-of-freedom (DoF) and corresponding shape functions at the elements cut by the material interface. Consequently, the elements cut by the interface are subdivided into sub-domains, each one of which is characterized by homogeneous material properties. Relying on the new nodal unknowns, G/X-FEM virtually introduces an interface-conforming mesh. Initially, this methodology was introduced for structural mechanics~\cite{duarte2000generalized,dolbow2000discontinuous,sukumar2003extended}, but later it was extended to analyze a wide range of problems in science and engineering~\cite{fries2010extended}. In ~\cite{aragon2010generalized,cervera2021comparative} the authors reported that the  G/XFEM method yields accuracy and convergence rates similar to those achieved by conforming meshes using standard FEM~\cite{aragon2010generalized,cervera2021comparative}. In the context of electro-statics/dynamics, Rochus et al.~\cite{rochus2011electrostatic} introduced the G/X-FEM for modeling interfaces between conductors and dielectric materials. Later, Soghrati et al.~\cite{soghrati2012interface} introduced a simplified G/X-FEM to handle discontinuous gradient fields. 

A drawback of G/X-FEM stems from the alteration of the size and graph of the system matrices due to introduction of additional degrees of freedom at the cut elements~\cite{zhang2016interface}. While for static problems this does not manifest any efficiency bottleneck as the cut elements do not alter, for cases where material interface evolves, the cut elements and, consequently, the additional degrees of freedom frequently change. 
Details regarding other pros and cons of G/XFEM can be found in the works of Soghrati et al.~\cite{soghrati2016automated} and Aragón et al.~\cite{aragon2020stability}.

An alternative to G/X-FEM is the Enriched Finite Element Method (E-FEM)~\cite{idelsohn2017elemental}, which shares the basis of the X-FEM method in terms of introducing additional degrees of freedom for representing the discontinuity. However, in contrast with G/X-FEM, these degrees of freedom are locally condensed (i.e. element-by-element) prior to the assembly of the global discrete system. This condensation becomes possible due to the appropriate choice of the enrichment shape functions, that must ensure independence among neighboring elements. This strongly improves the efficiency of the computations since enrichment does not affect the graph of the global system. This feature catapults the E-FEM as one of the most efficient techniques for the approximation of non-smooth solutions and is often used in the field of hydrodynamics for multi-fluid systems, where due to density differences discontinuous pressure gradients arise. Early contribution of the E-FEM for modeling discontinuous pressure field in multi-fluid problems was made in~\cite{coppola2005improving}. Later contribution to develop E-FEM in multi-phase flow  and thermal problems can find in \cite{ausas2010improved,hashemi2020enriched,hashemi2021three} and \cite{marti2017improved, tanyildiz2020solution}, respectively.
\\
Although E-FEM has been already used for several applications, including multi-phase flows and thermal applications, its employment to electromagnetic applications has scarcely been addressed in the literature to the best of our knowledge. In this sense, the present work makes the first step towards establishing a enriched finite element framework for such problems, aiming at accurate capturing of the discontinuity in the electric field in multi-material problems. 

In this work, we propose an easy-to-implement enrichment shape function capable of capturing discontinuity in the gradient of the field variable (electric potential), while not imposing any strong discontinuity. Using exclusively one enrichment degree of freedom per cut element is proposed. For the sake of simplicity the paper concentrates on electrostatic problems. Its generalization to other types of problems (dynamics, and electromagnetism) is straight-forward. The remaining part of the paper proceeds as follows: the basic governing equations of the model problem (multi-material electrostatics) description is presented in section 2. The enriched finite element formulation is presented in section 3. The numerical benchmarks are shown in section 4. Finally, summary and conclusions are presented.

\section{Governing equations}

\subsection{Problem statement}\label{sec_statement}
\begin{figure}
    \centering
                \includegraphics{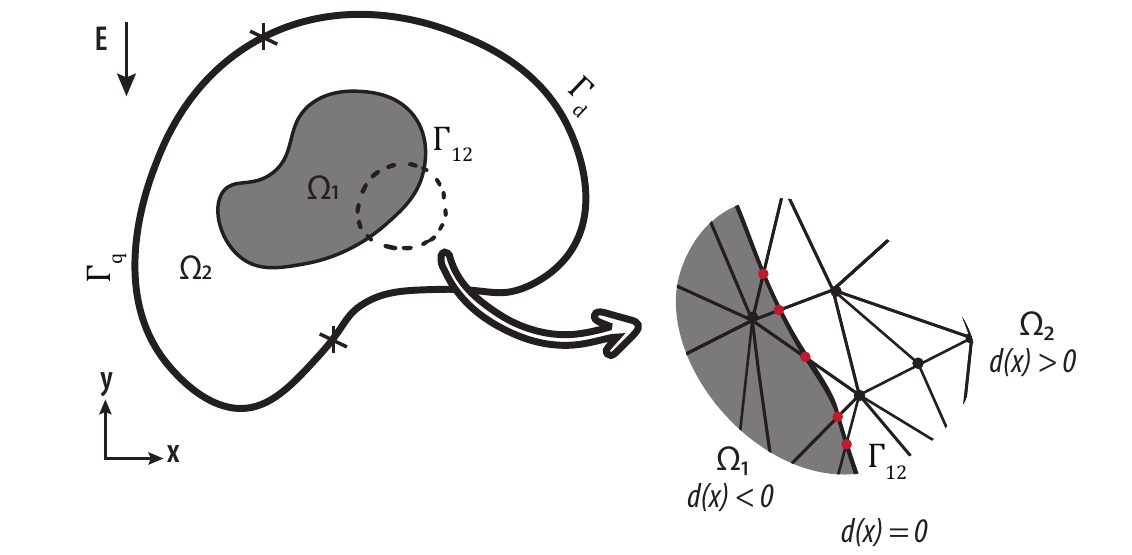}
                \caption{Schematic configuration of the geometry and boundary conditions used in-homogeneous electrostatic problems. The domain $\Omega$ is divided into two regions $\Omega_1$ and $\Omega_2$. \textit{Left} schematic description of the interface ($\Gamma_{12}$)} 
                \label{figure:domain}
\end{figure}
Let us consider a multi-material domain exposed to electric field (\textbf{E}) as shown schematically in Fig.~\ref{figure:domain}. This domain ($\Omega$) is composed of sub-domains $\Omega_1$ and $\Omega_2$, which represent two different media, each one with homogeneous properties. The governing equation of the electrostatic problem can be expressed as follows
\begin{equation}
    \nabla\cdot(\epsilon\textbf{E}) = 0
    \label{equation:weak_laplace}
\end{equation}
where $\epsilon$ represents the electrical material property specifying permittivity in case of dielectrics.
The electric field and the permittivity are related by the electric displacement \textbf{D} as ($\textbf{D} = \epsilon \textbf{E}$), which is the electric flux density between the media.

As an auxiliary variable that facilitates the definition of the problem, electric potential $\phi$ can be defined according to
\begin{equation}
    \textbf{E} = \nabla\phi.
    \label{equation:electricfield}
\end{equation}

\subsection{Discontinuity conditions in the electric potential gradient}
In the model problem described in Section \ref{sec_statement}, electric field \textbf{E} exhibits a discontinuity across the material interface due to an abrupt change in the material properties. Such discontinuity is shown in Fig.~\ref{figure:domain}.
In the present work we shall consider that the material properties are constant in each sub-domain and the interface is captured in a sharp manner (i.e. without smearing). We will therefore consider that the position of interface is defined by a smooth distance function $d(x)$ processing zero value at the interface. In order to distinguish between the two materials, a signed value of the distance function will be used: positive for the material on one side and negative for the material on the other side of the interface. 

\begin{figure}
    \centering
                \includegraphics[width=\textwidth]{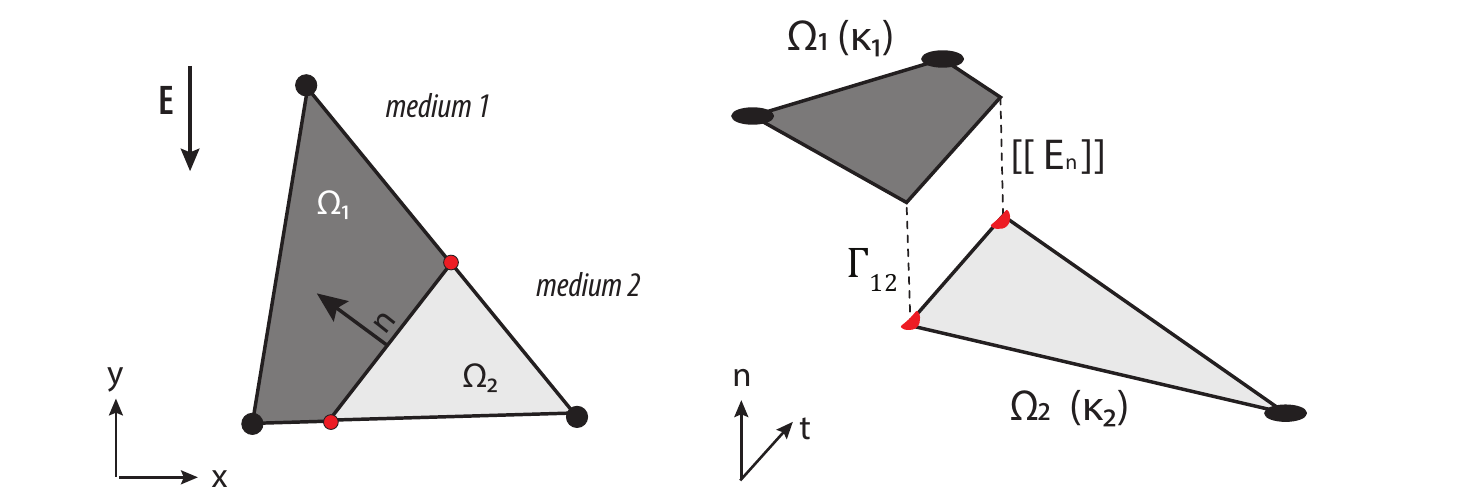}
                \caption{Schematic representation of the electric field at the cut elements. The standard nodes are represented with $\bullet$ symbol, while virtual nodes at the interface are represented with $\bullet$ symbol.} 
                \label{figure:electricfield_normal}
\end{figure}

The relationship between electric potential gradient (and, thus, the electric field) between both media is constrained by the ``jump condition'' ($[[\cdot]]$). Fig.~\ref{figure:electricfield_normal} represents this discontinuity along the interface, where \textbf{n} and \textbf{t} are the normal and tangential unit vector at the surface, respectively. It is worth mentioning that the electric potential itself, and consequently, the tangential component of the electric field $E_t = \textbf{E}\cdot\textbf{t}$ are continuous across the interface ($[[\phi]] = 0$ and $[[E_t]] = 0$), while the normal electric field $E_n = \textbf{E}\cdot\textbf{n}$ is discontinuous. 

\section{Finite element formulation}

The strong form of the equations governing multi-material electrostatics problems can be expressed as:
\begin{equation}    
\begin{matrix}
\nabla\cdot(\epsilon\nabla\phi) = 0  \quad on \quad \Omega_i \\
\phi = \widetilde{\phi}              \quad on \quad \Gamma_d \\
\mathbf{n} \cdot (\epsilon \nabla \phi) = 0  \quad on \quad \Gamma_{q}
\label{equation:strong_electricfield}   
\end{matrix}
\end{equation}
where $\widetilde{\phi}:\Gamma_d\rightarrow \mathbb{R}$ is the prescribed electric potential on the Dirichlet boundary condition while zero normal electric displacement is considered on Neumann boundary $\Gamma_q$. Given the test function $\omega \subset H^1(\Omega) = \{\omega : \phi|_{\Gamma_d} = 0 \}$), the weak form of Eq.~\ref{equation:strong_electricfield} can be written as 

\begin{equation}
\int_{\Omega} \nabla \omega \cdot\left(\epsilon\nabla \phi \right) d\Omega =  0
\label{equation:weakform} 
\end{equation}

In what follows, $\Omega$ is discretized into elements $\Omega_e$, such that ${\Omega} = \cup_e {\Omega_e}$ with $e = 1, \dots, n_{el}$ being the indices of the elements. Therefore, the integral over entire domain $\Omega$ that appears in Eq.~\ref{equation:weakform} can be expanded into the sum of elemental integrals as

\begin{equation}
\begin{matrix}
\sum^{n_{el}}_{e=1} \int_{\Omega_e} \nabla \omega \cdot\left(\epsilon\nabla \phi \right) d\Omega - \sum^{n_{el}}_{e=1}  \int_{\Gamma_{l}\cap\Omega_e} \omega \mathbf{n} \cdot \textbf{D}^+ d\Gamma = 0.
\end{matrix}
\label{equation:sum_weakform}
\end{equation}

Here, $\mathbf{n}$ is the outward normal to $\partial \Omega_e$. The second term in Eq.~\ref{equation:sum_weakform} incorporates the contribution of the electric displacements from neighboring elements ($\textbf{D}^+$) along the inter-elemental boundaries ($\Gamma_{l}$). In other words, this term guarantees the adaptation of the electric potential within the element to the inter-elemental Neumann condition arising from the adjacent cut elements. 
Using continuous polynomial approximation for the variables, as in the standard FEM, the contribution of the inter-elemental normal electric displacements is counter-balanced during the assembly procedure, making this term vanish in the assembled system.
However, in case the continuity is disturbed, the incorporation of this term is necessary.
In the present work, as it will be shown in Section 3.2, this appears to be the case in the cut elements, making accounting for the second term in Eq.~\ref{equation:sum_weakform} necessary.

\subsection{Enriched finite element space}

Although applying standard finite element approximations (continuous interpolations for the test and trial functions) to the Eq.~\ref{equation:sum_weakform} can be used to approximate the solution of an electrostatic problem in a multi-material domain, it would lead to important inaccuracies associated with the cut elements. The main problem originates from the fact that the mentioned equation can represent continuous gradients only, while in the problem at hand gradient field passing across the material interface is no longer continuous \cite{soghrati2012interface}. Thus, standard finite element method applied to  Eq.~\ref{equation:sum_weakform} would imply approximating discontinuous gradient by a continuous one, which is unacceptable for the situations where the material properties on the two sides of the interface strongly differ. 

This issue can be addressed by enriching the solution space. In the present work we propose to enrich the finite element space by adding exclusively one degree of freedom, responsible for controlling the size of the ``jump'' in the gradient of the potential across the cut. This setting would allow to subsequently condense the enrichment contribution, being completely local to the element.

The proposed approximation to the solution can thus be expressed at the level of a single cut element as:

\begin{equation}
\phi^h(x) = \underbrace{\sum_{i\in e.n} N_i(x) \phi_i}_{standard} +
\underbrace{ \Bar{N}(x) \phi^*}_{enriched}
    \label{equation:enriched_approx_func}
\end{equation}

\begin{equation}
    \omega^h(x) = \sum_{i\in e.n} N_i(x) \omega_i + \Bar{N}(x) \omega^*
    \label{equation:enriched_weight_func}
\end{equation}

The first terms in these equations represent the standard FEM contribution with $N_i$ being the standard shape function and $\phi_i$ and $\omega_i$ the values corresponding to node $i$ of the element. The enrichment terms are defined as the product of the enrichment shape function $\Bar{N}(x)$ and the additional degrees of freedom, $\phi^*$ and $\omega^*$, controlling the jump in the gradient at the cut. 

Taking into account its continuity property, we propose to rely on the hat-function for the definition of the enrichment shape function as
\begin{equation}
    \Bar{N}(x) = \sum_{i\in e.n} N_i(x) \left|d_i(x)\right| - |\sum_{i\in e.n} N_i(x) d_i(x)|
\label{equation:enrichmentfunction}
\end{equation}
\begin{figure}
\centering
\includegraphics[width=\textwidth]{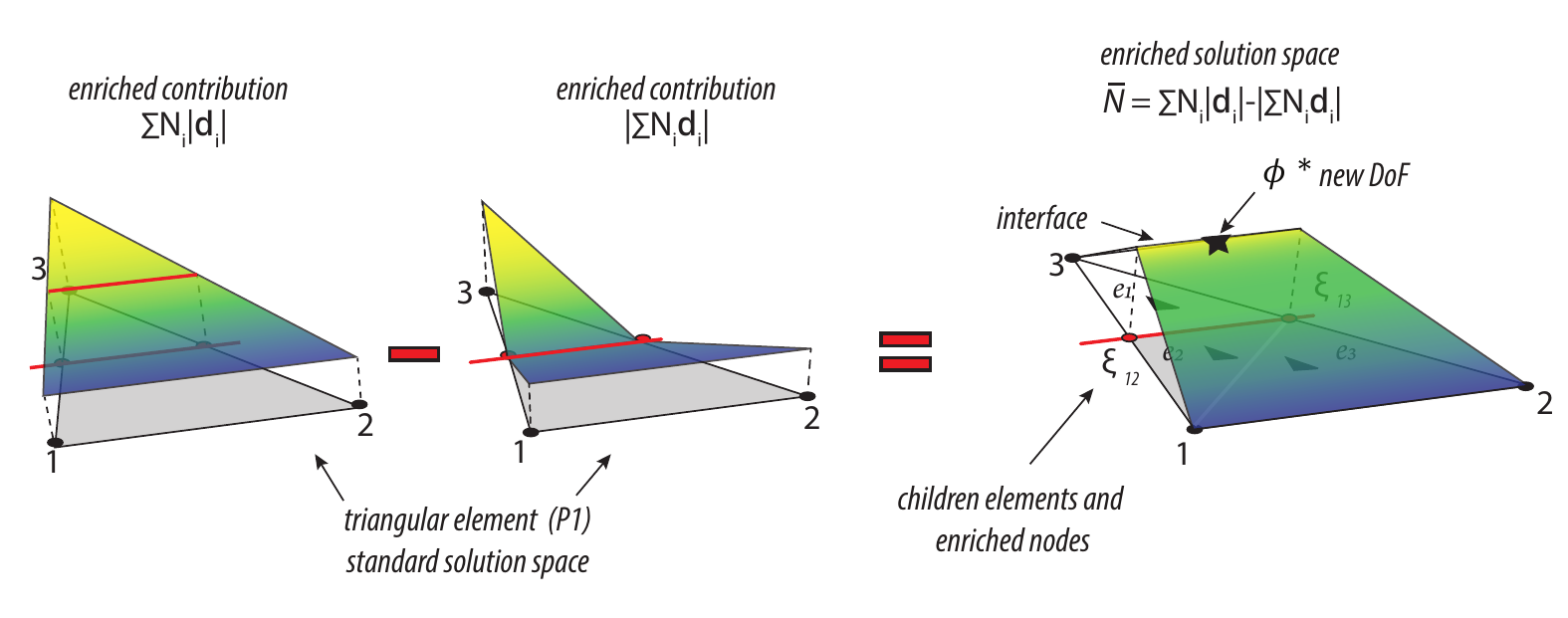}
\caption{Construction of the enriched space on P1 triangular element. \textit{Left} schematic representation of each term of the hat-function. \textit{Right} enriched finite element, the parent element is split into children elements ($e_1$, $e_2,$ and $e_3$) and virtual nodes at the interface ($\xi$). The Gauss points of the children elements and the additional DoF are represented with $\blacktriangle$ and $\star$ symbols, respectively.}
\label{figure:split_element}
\end{figure}
\\
Figure~\ref{figure:split_element} shows a schematic of the enrichment shape function (Eq.~\ref{equation:enrichmentfunction}). It can be seen that the function is weakly discontinuous across the interface and is smooth elsewhere, which guarantees that all nodal values of the background solution space are zero. In the present study P1 triangular elements are used, which are subdivided into "children" elements, as shown in Fig.~\ref{figure:split_element}. 

The hat-function defined in Eq. \ref{equation:enrichmentfunction} was first used by M\"oes et al.~\cite{moes2003computational} in the context of G/X-FEM for modeling complex  microstructures in elasticity problems. It has been later applied for the G/X-FEM modeling of various bimaterial interface problems~\cite{fries2010extended}. It was shown that utilizing this approach, optimal convergence rates can be obtained for linear and bilinear elements.
However, in the context of the G/X-FEM, in order to produce the enriched shape functions, the hat-function is usually multiplied by the corresponding standard nodal shape functions in order to preserve the compatibility between the neighboring elements.
Such shape functions require the introduction of the new nodal degrees of freedom.
Therefore, being incapable of the static condensation of additional degrees of freedoms for the elemental system of equation, the development of an efficient E-FEM framework is hindered.
In the present work, as already mentioned, we propose to use a single (elemental) degree of freedom facilitating the static condensation. On the other hand, the compatibility between the neighboring cut elements will be enforced without introducing additional degrees of freedom.

\subsection{Weak formulation of the enriched space}
\label{section_enriched_formulation}

The residual form of the enriched finite element approximation of the governing equations can be obtained by substituting the approximation of the test and trial functions (Eqs.~(\ref{equation:enriched_approx_func}) and~(\ref{equation:enriched_weight_func})) into the weak form defined in Eq. ~(\ref{equation:sum_weakform}).

\begin{equation}
\begin{aligned}
\mathcal{R}(\phi^h, \omega^h) & = \sum^{n_{el}}_{e=1} \int_{\Omega_e} \nabla \left( \sum_{i\in e.n} N_i(x) \omega_i + \Bar{N}(x) \omega^*\right)\cdot[\epsilon\nabla( \sum_{j\in e.n} N_j(x) \phi_j + 
\\
& \Bar{N}(x) \phi^*)] d\Omega  - \sum^{n_{el}}_{e=1}  \int_{\Gamma_{l}\cap\Omega_e} \Bar{N}(x) \omega^* \mathbf{n} \cdot \textbf{D}^+ d\Gamma = 0.
\end{aligned}
\label{equation:sum_weakform_residual} 
\end{equation}

Note that the continuity of the electric displacement ($\textbf{D}^+$) leads to counter-balance of this term between the neighboring elements across the shared edges. Upon assembly of the elemental contributions this term should vanish and therefore, is not typically included in the implementation of the elemental equations. However, the test function corresponding to the enriched space ($\bar{N}(x)\omega^*$) proposed in this work is local to the element. Therefore, its contribution to the assembled system of equation cannot be simply omitted. 

\begin{figure}[h]
    \centering
                \includegraphics[width=\textwidth]{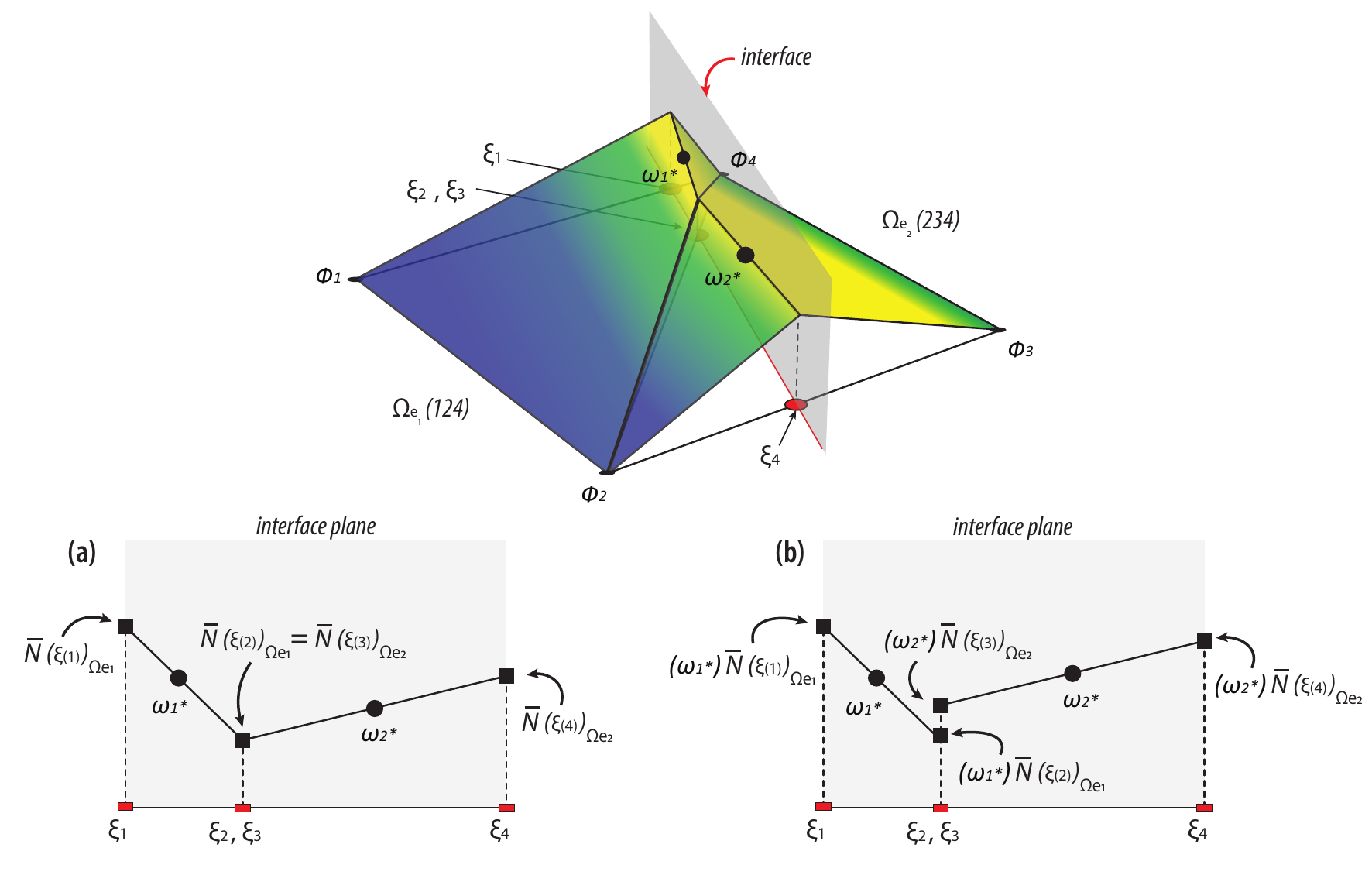}
                \caption{Simplified representation of enriched space between neighbor elements cut by an orthogonal plane to the interface. (a) enrichment function. (b) product of the enrichment function and added DoF in the solution space. The $\xi$ represents the virtual nodes at the interface} 
                \label{figure:discontinuous_enrich}
\end{figure}

To have a closer look at the above-mentioned issue originating from the local nature of the enrichment, let us consider two neighboring elements cut by a plane orthogonal to the interface (see Fig.~\ref{figure:discontinuous_enrich} [a]). One can see that even though the enrichment shape functions themselves satisfy the compatibility condition everywhere, this condition on the cut edges may not me guaranteed for the enriched field approximation, which is local to each cut element (see. Fig.~\ref{figure:discontinuous_enrich} [b]). 

The last term in Eq. \ref{equation:sum_weakform_residual} accounts for the continuity of the electric displacement, which is equivalent to imposing the so-called \textit{inter-element dissipate flux continuity} in diffusion problems. Omitting this term in the context of enriched finite element approaches compromises the accuracy of the method as shown in \cite{idelsohn2017elemental} and \cite{marti2017improved}.  

The matrix form of the governing equations is obtained by standard procedure computing the derivatives of  the residual (Eq.~\ref{equation:sum_weakform_residual}) with respect to $\omega_i$ and $\omega^*$ and, subsequently, with respect to $\phi_i$ and $\phi^*$. The discrete system of equations can be thus written as:
\\
\begin{equation}
\begin{bmatrix} 
\textbf{K} & \textbf{B} \\
\textbf{B}^T  - \textbf{D} & \textbf{K}_{enr} - \textbf{D}_{enr}
\end{bmatrix}
\begin{bmatrix} 
\boldsymbol{\phi} \\ \phi^* 
\end{bmatrix}
=
\begin{bmatrix} 
0 \\ 0
\end{bmatrix}
\label{equation:system_eq}
\end{equation}
where matrix $\textbf{K}$, vectors $\textbf{B}$ and $\textbf{D}$, and scalars $\textbf{K}_{enr}$ and $\textbf{D}_{enr}$ are defined as follows

\begin{equation*}
\begin{bmatrix}
    \textbf{K}_{ij}= \int_{\Omega_e} \nabla N_i \cdot\left(\epsilon\nabla N_j \right) d\Omega  \qquad  &   \textbf{B}_{i}= \int_{\Omega_e} \nabla N_i \cdot\left(\epsilon\nabla \Bar{N} \right) d\Omega
\end{bmatrix}
\end{equation*}
\begin{equation*}
\begin{bmatrix}
    {K}_{enr}= \int_{\Omega_e} \nabla \Bar{N} \cdot\left(\epsilon\nabla \Bar{N} \right) d\Omega  
\end{bmatrix}
\end{equation*}
\begin{equation*}
\begin{bmatrix}
    \textbf{D}_{i}= \int_{\Gamma_{l}\cap\Omega_e} \Bar{N} \textbf{n} \cdot \left(\epsilon \nabla N_i\right) d\Gamma  \qquad  &  \textbf{D}_{enr}= \int_{\Gamma_{l}\cap\Omega_e} \Bar{N} \textbf{n} \cdot \left(\epsilon \nabla \Bar{N}\right)d\Gamma
\end{bmatrix}
\label{equation:discerete_final}
\end{equation*}

In the implementation of the present work, the enriched DoF $\phi^*$ associated to the weak discontinuity of the electric potential is condensed at the elemental level

\begin{equation}
    \left(\textbf{K} - \textbf{B}[\textbf{K}_{enr}-\textbf{D}]^{-1}[\textbf{B}^T-D] \right)\boldsymbol{\phi} = 0,
    \label{equation:phi_condensed}
\end{equation}

so that the discontinuity is captured, without affecting size and graph of the global system.

After solving the complete system defined by Eq.~\ref{equation:system_eq}, $\phi^*$ can be recovered by
\begin{equation}
    \phi^{*} = -[\textbf{K}_{enr}-\textbf{D}]^{-1}[\textbf{B}^T-D] \boldsymbol{\phi}
    \label{equation:phi_enrich_decondensed}
\end{equation}

Note that recovery of the enrichment degree of freedom after the obtaining the solution is necessary only for the calculation of the electric field at the post-processing step.

\textcolor{Blue}{\section{Numerical examples}}
In this section the performance of the proposed enriched finite element method is tested. Three examples are provided to assess the performance of the method, paying particular attention to its accuracy. The benefits of accounting for the electric displacement term are also assessed. 

The method is implemented in \textit{Kratos Multiphysics}, an in-house Open Source C++ object-oriented Finite Element framework~\cite{kratosweb, dadvand2010object}. For the solution of the linear system a stabilized Bi-conjugate gradient solver (BiCGSTAB) was applied \cite{sleijpen1993bicgstab}. The convergence tolerance was set to  $10^{-8}$.

\subsection{Two materials between horizontal electrodes in a square domain}
\begin{figure}[h]
    \centering
                \includegraphics[width=\textwidth]{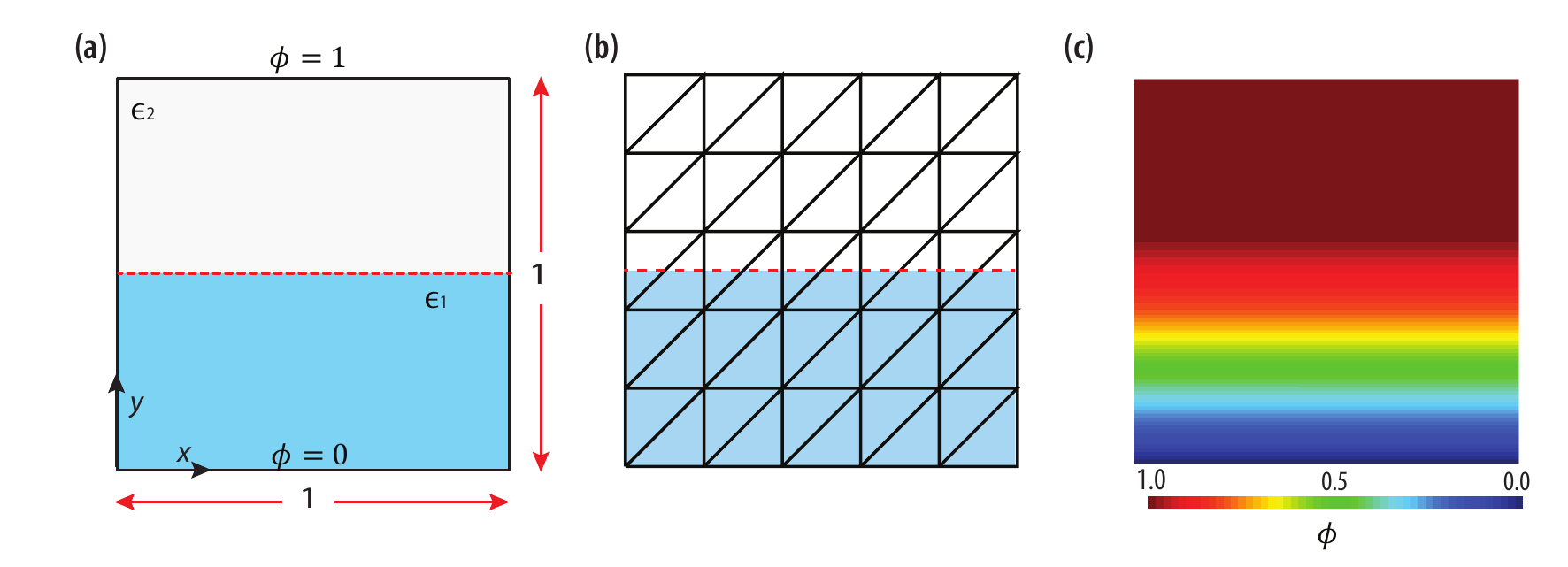}
                \caption{The planar test for bi-material problem. (a) Sketch of the problem composed of two materials with electric properties $\epsilon_1$ and $\epsilon_2$. (b) Nonconforming mesh with element size \textit{h} = 0.2,  employed for standard FEM and E-FEM. (c) Approximation of the electric potential field for $1/Q\approx0$ using E-FEM with the inclusion of inter-elemental D$^+$ term in the tangent matrix.} 
                \label{figure:1d_domain}
\end{figure}
The configuration used in the first test is shown in Fig.~\ref{figure:1d_domain}: the electric potentials of 1 $V$ and 0 $V$ are applied to two parallel electrodes coinciding with the top and bottom walls of a unit square domain, respectively. The domain contains two homogeneous materials, on both sides of the horizontal interface located at the mid-height of the domain. 

The numerical simulation was performed for three different ratios of the material properties (permittivities)  $Q = \epsilon_1 / \epsilon_2 $ corresponding to the upper and lower sub-domains: $Q = 1$, $Q = 3$, and $1/Q \approx 0$, which corresponds to the \textit{conducting-dielectric} case-study. Intentionally, a very coarse mesh was used in the simulation (see Fig.~\ref{figure:1d_domain} [b]).

The problem was solved 
\begin{itemize}
    \item  in the absence of enrichment function (standard FEM)
    \item  with enrichment (but without the electric displacement term (E-FEM no D$^+$))
    \item  with enrichment and the electric displacement term (E-FEM incl. D$^+$)
\end{itemize}

\begin{figure}
    \centering
                \includegraphics[width=\textwidth]{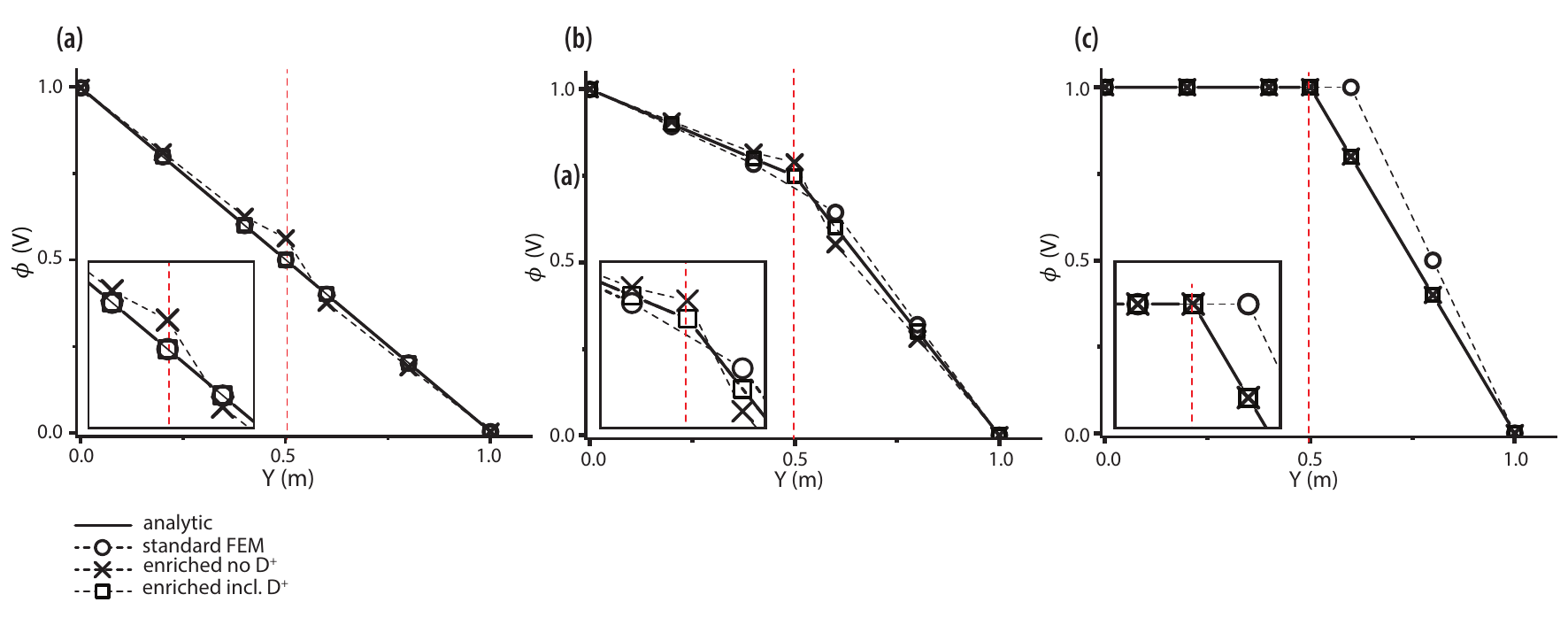}
                \caption{Distribution of the electric potential along the vertical cut. For three different ratios of electric properties: (a) $Q = 1$, (b) $Q = 3$, and (c) $1/Q \approx 0$.} 
                \label{figure:validation_y_dir}
\end{figure}

\textit{Electric potential approximation}\\
Fig.~\ref{figure:validation_y_dir} shows the comparison between the results obtained using the standard FEM, the proposed E-FEM and the analytic solution (see Table~\ref{table:phi_sol}). Electric potential $\phi$ is plotted along the vertical cut made in the middle of the domain (x=0.5 m). Note that in order to visualize the results of the E-FEM, the additional degree of freedom was "recovered" after solving the system, following Eq. (\ref{equation:phi_enrich_decondensed}). This way, the value exactly at the cut (where no actual node is located) could be reconstructed as shown in the inset of the Fig.~\ref{figure:validation_y_dir}. 
\begingroup
\begin{table}
  \centering 
\caption {\label{table:phi_sol} Analytic dimensionless solution for the bi-material planar cases. The subscript $1$ and $2$ correspond to medium 1 and 2, respectively. $\textbf{E}_\infty$ is equal to -1.}
\begin{tabular}{cllll}
\hline
Case & $\phi_1$ & $\phi_2$ & $E_1$ & $E_2$ \\
\midrule
\textit{Q = 1} & $1-y$ & $1-y$ & $E_\infty$ & $E_\infty$\\
\textit{Q = 3} & $3y/2$ & $(y+1)/2$ & $1.5 E_\infty$ & $0.5 E_\infty$\\
\textit{1/Q $\approx$ 0} & $2y$ & 1 & $2E_\infty$ & 0\\
\midrule
\end{tabular}
\end{table}
\endgroup
\\
In this test case (which is essentially uni-dimensional), according to the exact solution (Table~\ref{table:phi_sol}) electric potential must decrease linearly along the Y direction. The slope changes precisely at the cut (proportional to the permittivity value), and is constant on each side of the interface.
\\
One can see that the standard FEM fits the exact solution only for the single material case (i.e. when both sub-domains have the same material properties ($Q = 1$), see Fig.~\ref{figure:validation_y_dir} [a]). However, for ratios greater than 1 ($Q > 1$) the standard FEM fails to capture the electric potential at the interface. The inaccuracy becomes even more evident in the case when one of the sub-domains contains a conductor media  ($1/Q \approx 0$) (see Fig.~\ref{figure:validation_y_dir} [c]). Although the value of $\phi$ may be close to the correct one, standard FEM is incapable of accounting for the bi-material condition in the cut elements. 
\\
As expected the proposed E-FEM improves the results for bi-material cases. However, the results put in evidence the importance of accounting for the inter-elemental electric displacement term, especially in the case of small ratios. If this term is not included in the tangent matrix, the results accuracy becomes compromised.
This is particularly evident in Fig.~\ref{figure:validation_y_dir} [a], for the limiting case of $Q = 1$. Similar results were obtained when the ratio increased to 3 ($Q = 3$), both standard FEM and E-FEM without the inclusion of D$^+$ term fail to approximate the potential (see Fig.~\ref{figure:validation_y_dir} [b]). From these data, it is apparent that accounting for the inter-elemental electric displacement condition between adjacent cut elements is essential for obtaining accurate approximation of the electric potential at the interface in the proposed E-FEM framework. 
\begin{figure}
    \centering
                \includegraphics[width=\textwidth]{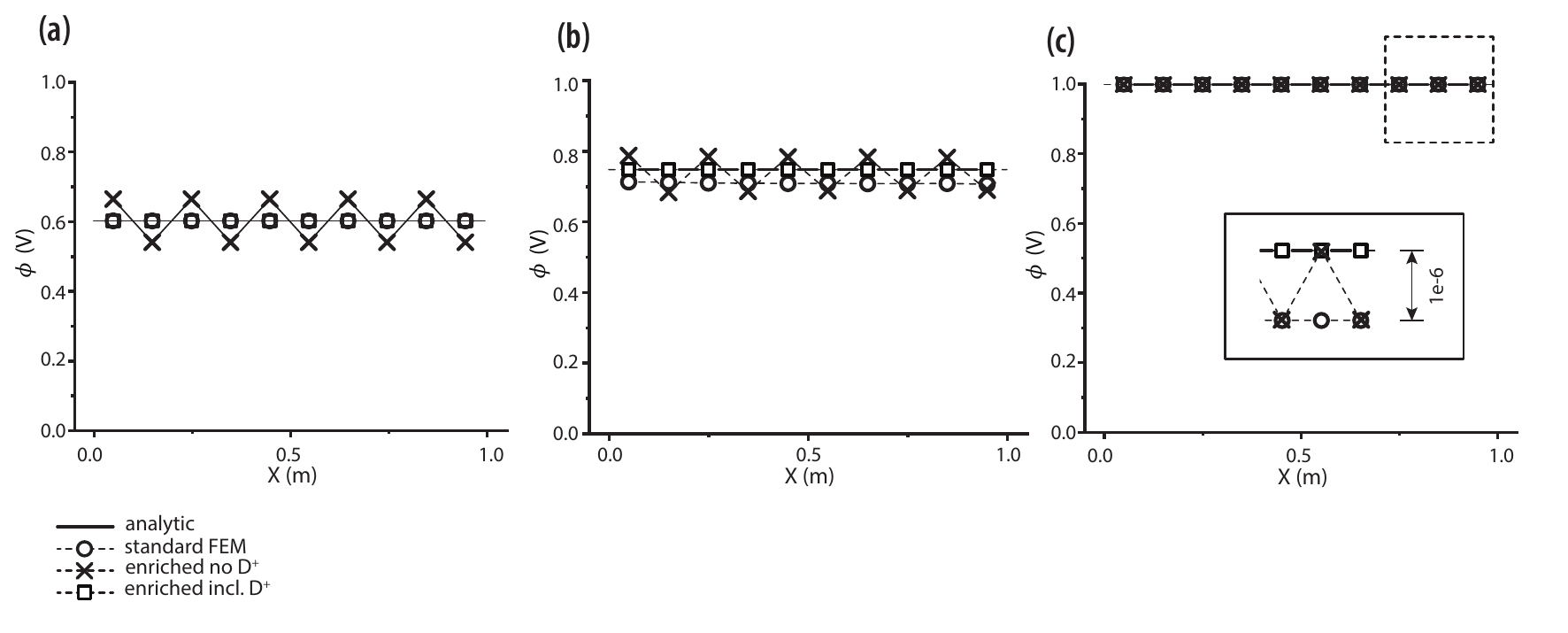}
                \caption{Distribution of the electric potential along the interface (Y=0.5). For three different ratios of electric properties: (a) $Q = 1$, (b) $Q = 3$, and (c) $1/Q \approx 0$.} 
                \label{figure:validation_x_dir}
\end{figure}

In order to have a closer look at the above-described phenomenon, the potential across the horizontal cut is displayed next. Figure~\ref{figure:validation_x_dir} [a] shows that the mentioned non-conformity affects the solution even when the same material is considered for both sub-domains. As the magnitude of $Q$ increases, the potential fluctuation diminishes, reaching a nearly-zero value for the case of conducting-dielectric media ($1/Q\approx0$). For instance, the difference between the potential at the interface of the neighboring elements decreased to $\sim 1e^{-6}$ when the conducting-dielectric media was considered (see Fig.~\ref{figure:validation_x_dir} [c]). What can be clearly seen from the graphs above is that the spurious variability of the potential is suppressed via including the inter-elemental electric displacement term into the enriched formulation. 
Quantitative results are provided in Table~\ref{table:phi_error}.

\begingroup
\begin{table}
  \centering 
\caption {\label{table:phi_error} Deviation from the exact solution of the potential $\phi$ at the interface: standard FEM and E-FEM with/without the inter-elemental electric displace term.}
\begin{tabular}{clll}
\hline
\multicolumn{1}{c}{element size}   & \multicolumn{3}{c}{error (\%) $|1-\phi_{interface}/\phi_{interface}^{ex}|$  } \\
\cline{2-4} 
($h$)   & standard FEM   &   E-FEM no D$^+$      &  E-FEM incl. D$^+$           \\
\midrule
\textit{Q = 1}\\
0.3 & 0 & 0.2 & 0 \\
0.2 & 0 & 0.125 & 0 \\
0.03 & 0 & 0.02 & 0 \\
\textit{Q = 3}\\
0.3 & 0.080 & 0.083 & 0 \\
0.2 & 0.046 & 0.077 & 0 \\
0.03 & 0.007 & 0.007 & 0 \\
\midrule
\end{tabular}
\end{table}
\endgroup

\textit{Approximation of the electric field}\\
While analyzing the approximation we noted the considerable improvement in the accuracy of the method in case of using E-FEM, its benefits become even more evident when looking at the electric field, which is proportional to the gradient of the potential and is discontinuous in a multi-material problem, which cannot be reproduced by the standard FEM. In order to illustrate it, the same test as before (see Fig.~\ref{figure:1d_domain}) is solved, however this time electric field is analyzed.

\begin{figure}[h]
    \centering
                \includegraphics[width=\textwidth]{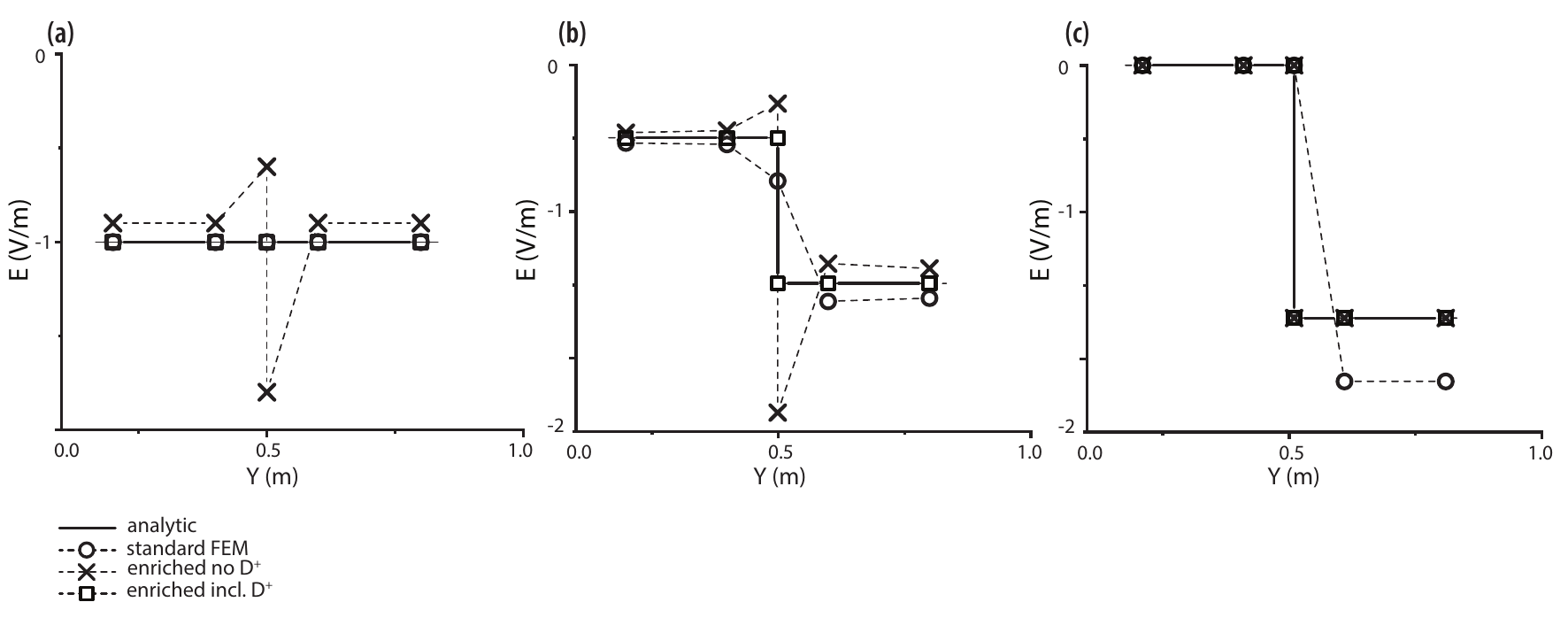}
                \caption{Distribution of the electric field along the vertical cut. For three different ratios of electric properties: (a) $Q = 1$, (b) $Q = 3$, and (c) $1/Q \approx 0$.} 
                \label{figure:validation_E}
\end{figure}

Figure~\ref{figure:validation_E} shows the electric field computed along the vertical cut made in the middle of the domain (x=0.5 m). The numerical results are plotted for three ratios of permittivities ($Q=1$, $Q=3$, and $1/Q\approx0$), adjoining its analytical solution that is presented in Table~\ref{table:phi_sol}. What stands out from Fig.~\ref{figure:validation_E} [a] is that for $Q = 1$ there should be no jump at the interface since both sides of the domain share the same electric property, however, the lack of the electric displacement term in the tangent matrix generates an artificial jump in the vicinity of the interface. 
In contrast, as expected, standard FEM solutions fail to capture the sharp jump of the electric field for $Q \neq 1$ (see Figs.~\ref{figure:validation_E} [b-c]). This outcome is a consequence of the continuity of the standard finite element shape functions, and, consequently, inability to represent inter-elemental discontinuities. \\

The inclusion of the electric displacement term in the tangent matrix enhances the numerical performance of E-FEM. One can see, that in this case E-FEM does not provide deviations from the exact solution (see Table~\ref{table:error_E}). Note that the errors of E-FEM without the inclusion of D$^+$ term show a dependence on the magnitude of the properties ratio: for very high ratios the electric displacement term may be limited or neglected. 

\begingroup
\begin{table}
  \centering 
\caption {\label{table:error_E} Deviations from the exact solution of \textbf{E} at the cut elements: standard FEM, E-FEM with/without the inter-elemental electric displace term.} 
\begin{tabular}{c cccccc}
\hline
\multicolumn{1}{c}{element size}   & \multicolumn{3}{c}{error (\%) $|1-E/E^{ex}|$  } \\
\cline{2-7} 
($h$)   & \multicolumn{2}{c}{standard FEM}    &   \multicolumn{2}{c}{E-FEM no D$^+$}     &  \multicolumn{2}{c}{E-FEM incl. D$^+$} \\
& $E_1$ & $E_2$ & $E_1$ & $E_2$ & $E_1$ & $E_2$ \\
\midrule
\textit{Q = 1}\\
0.3 & 0 & 0 & 0.8 & 0.4 & 0 & 0 \\
0.2 & 0 & 0 & 0.875 & 0.375 & 0 & 0 \\
0.03 & 0 & 0 & 0.98 & 0.34 & 0 & 0 \\
\textit{Q = 3}\\
0.3 & 0.512 & 0.463 & 0.522 & 0.493 & 0 & 0 \\
0.2 & 0.534 & 0.397 & 0.279 & 1.161 & 0 & 0 \\
0.03 & 0.477 & 0.570 & 0.441 & 0.519 & 0 & 0 \\
\textit{1/Q $\approx$ 0}\\
0.3 & 0.999 & - & 1e-5 & - & 0 & - \\
0.2 & 0.999 & - & 5e-6 & - & 0 & - \\
0.03 & 0.999 & - & 0 & - & 0 & - \\
\midrule
\end{tabular}
\end{table}
\endgroup
\subsection{Two-material problem with an inclined interface}
In the next example we analyze a truly two-dimensional case, studying a square domain containing two materials with non-horizontal interface. 
Besides the non-perpendicularity of the interface and the external electric field imposed by the Dirichlet boundaries, in this test, the topology of the individual material domains varies between the neighboring cut elements.
Here, the order of accuracy of the method is assessed for different permittivity ratios.

\begin{figure}
    \centering
                \includegraphics[width=12cm]{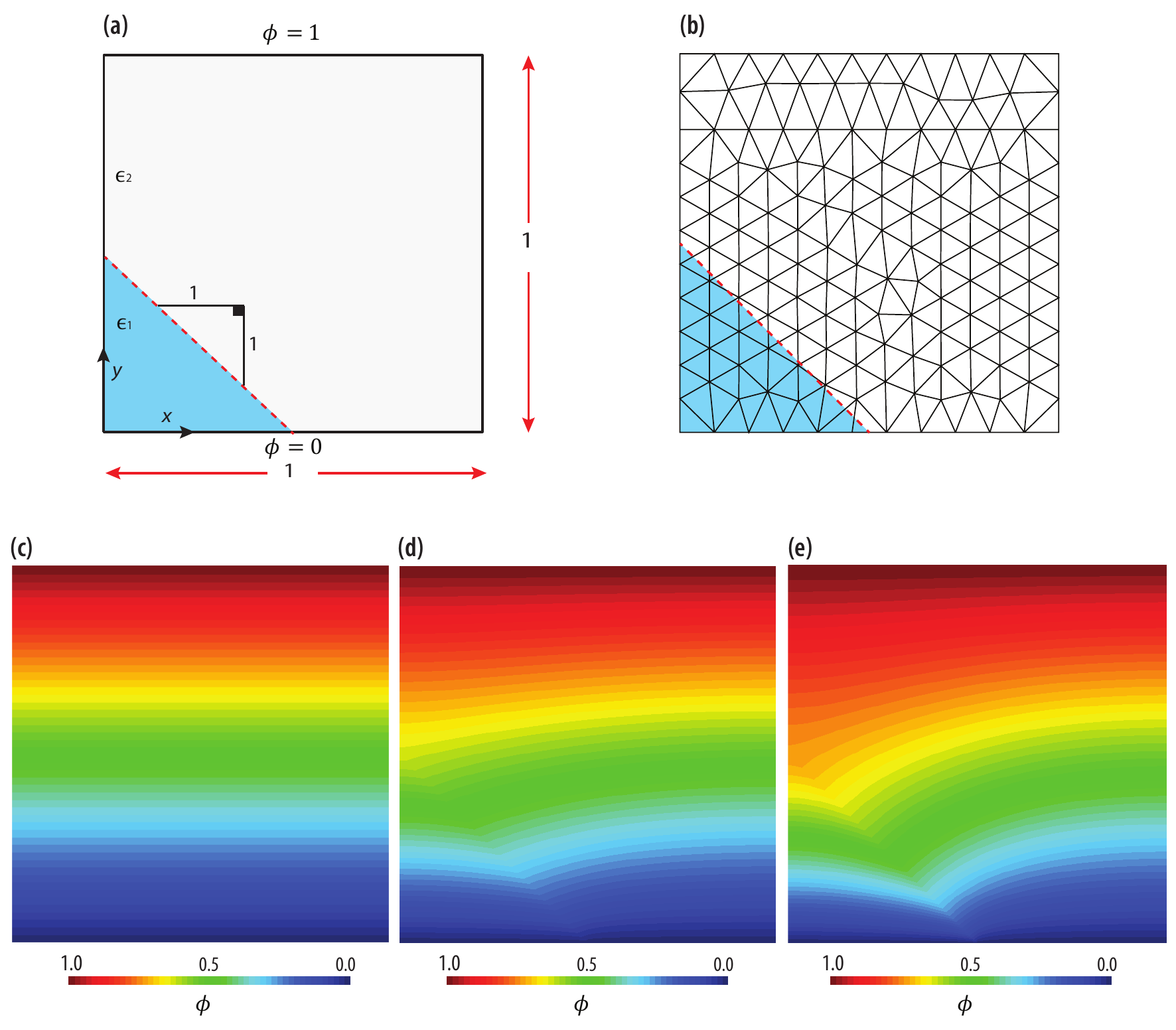}
                \caption{Model problem used to evaluate the electric displacement term in a case with inclined interface. a) Sketch of the problem composed of two materials with electric properties $\epsilon_1$ and $\epsilon_2$, and the interface slope of 1. b) Nonconforming mesh with element
size \textit{h} = 0.15. Electric potential distribution (reference solutions): c)  $Q = 1$, d) $Q = 3$, and $1/Q \approx 0$} 
                \label{figure:mesh_45_deg}
\end{figure}

The test setting is shown in Fig.~\ref{figure:mesh_45_deg}. The domain contains two materials separated by an inclined interface forming a slope of 45$^o$. The domain size and boundary conditions are identical to those used in the previous example: the electric potentials of $1V$ and $0V$ are applied to the top and bottom walls of the domain, respectively. 
The same set permittivity ratios is kept for this example ($Q=1 $, $Q=3 $, and $1/Q \approx 0 $). The problem is solved in the absence of the electric displacement term and with its contribution to the tangent matrix. Since no analytic solution exists for this example, the results obtained using standard FEM on a fine interface-conforming mesh with $h =0.01$ are used as the reference solution presented in Fig.~\ref{figure:mesh_45_deg}[c], \ref{figure:mesh_45_deg}[d], and \ref{figure:mesh_45_deg}[e], for different permittivity ratios.

\begin{figure}
    \centering
                \includegraphics[width=\textwidth]{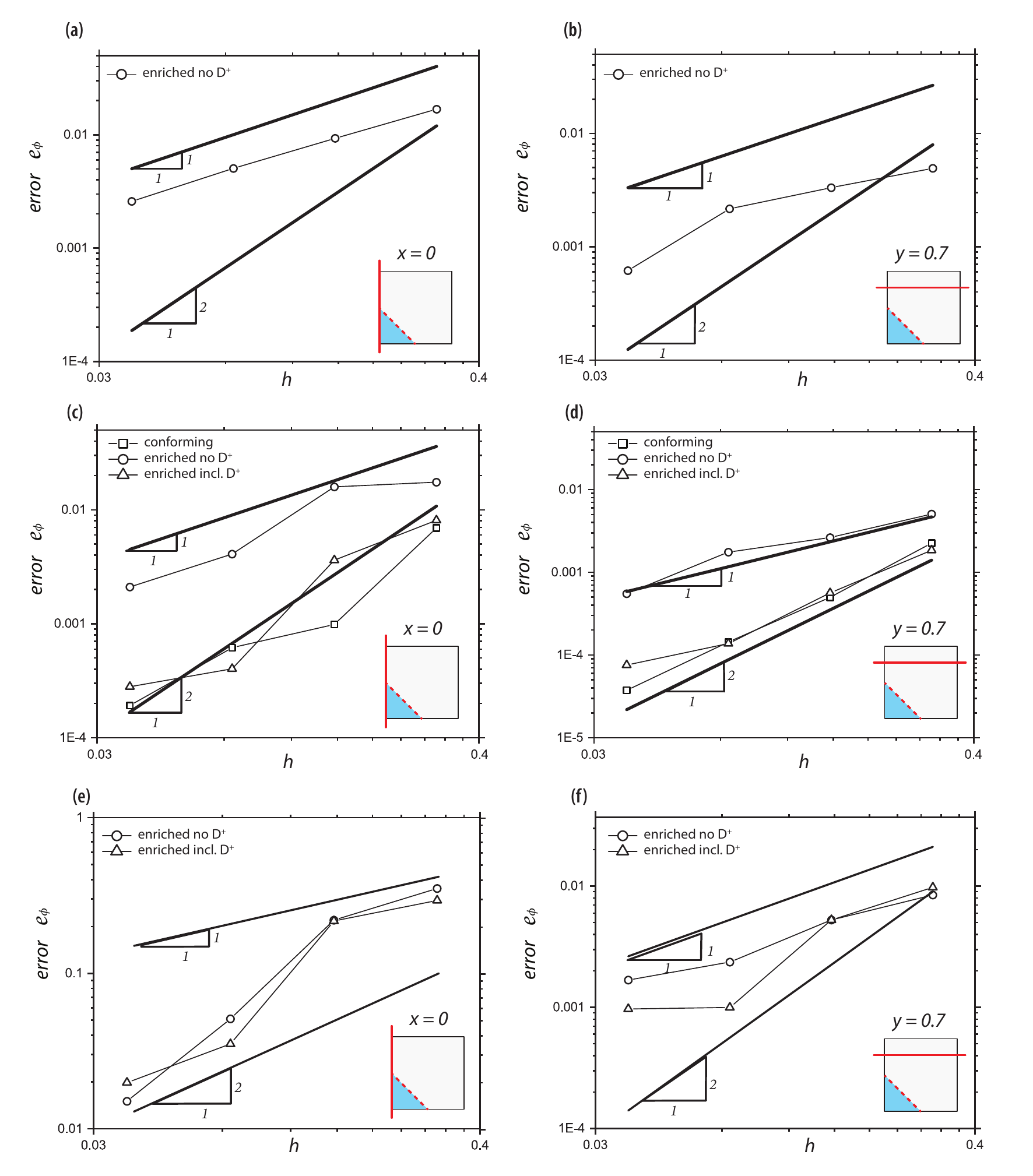}
                \caption{The rate of error in $\phi$ versus mesh-size along X = 0 and Y = 0.7 obtained for $Q = 1$ (a-b), $Q = 3$ (c-d), and $1/Q \approx 1000$ (e-f).} 
                \label{figure:error_plot.eps}
\end{figure}

In order to estimate the accuracy of the method, error is analyzed as a function of the element size. In this work, error (in $L_2$ norm) is measured as the root-mean-square (RMS) of the difference between the obtained and the reference solution integrated along a given line passing through the domain ($\Gamma$) as

\begin{equation}
    e_{\phi} = \sqrt{\int_\Gamma (\phi^h - \phi^{ref})^2d\Gamma}
    \label{equation:error_integral_eq}
\end{equation}

E-FEM results are obtained on meshes with characteristic element sizes of $h = 0.3 , 0.15, 0.075$, and $0.0375$. 
Fig.~\ref{figure:error_plot.eps} shows the error for different permittivity ratios estimated for two transect regions, i.e. the values of $\phi$ are calculated over the vertical and horizontal cut lines of $x=0$ and $y=0.7$. Note that in order to visualize the effect of the spurious values on the non-split elements, the horizontal cut was taken far from the interface region.

For the particular case of $Q = 1$, the E-FEM including the inter-elemental D$^+$ term yields the same result as the reference solution, thus, no deviations are observed in Fig.~\ref{figure:error_plot.eps} [a-b]. In contrast, E-FEM without the inclusion of the inter-elemental electric displacement term in the tangent matrix leads to a convergence rate of $1.0$, even in the non-split elements. One can see that by increasing the permittivity ratio, in Fig.~\ref{figure:error_plot.eps} [c-d], the convergence rate of E-FEM with the electric displacement term is $2.0$ that is similar to the one obtained with standard FEM on a conforming mesh (with the same characteristic element size), whereas the lack of the D$^+$ term deteriorates the performance of E-FEM generating a convergence rate of almost first-order. \\

Figure~\ref{figure:error_plot.eps} [e-f] shows the convergence rate for $1/Q\approx 0$, for which both enrichment methods (with/without inclusion of the inter-elemental electric displacement term) lead to similar convergence rates, which is slightly below second-order. These results confirm that the term of the electric displacement is limited by the magnitude of the permittivity ratio, or in other words, the benefits of including inter-elemental electric displacement term diminishes as $Q$ grows. Overall, one can see that the proposed E-FEM exhibits attractive convergence rates.

\subsection{Electric field around a cylinder}
\begin{figure}
    \centering
                \includegraphics[width=\textwidth]{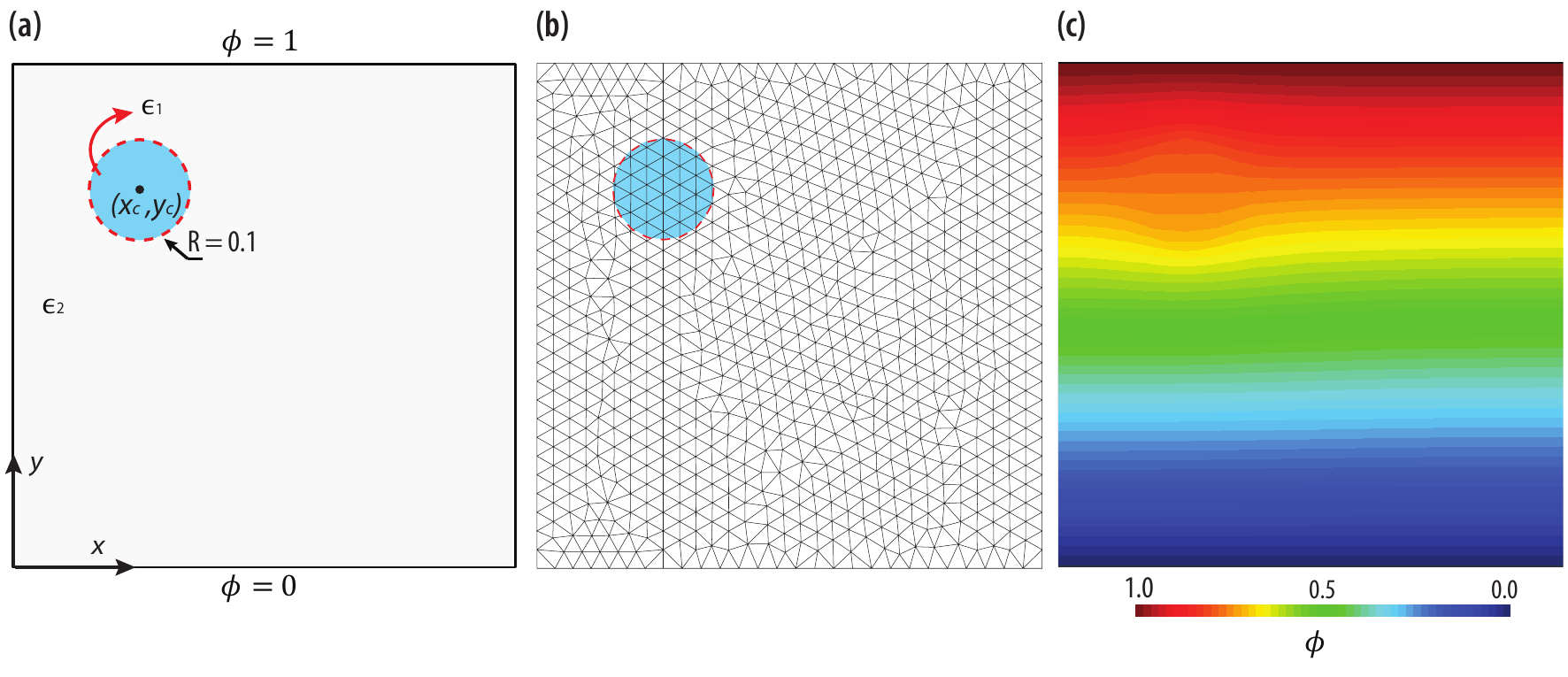}
                \caption{Model problem with a curved interface. a) Domain and boundary conditions. b) Nonconforming mesh with element
size\textit{ h} = 0.0375. c) Electric potential distribution obtained using the proposed E-FEM with the inclusion of D$^+$ term.} 
                \label{figure:circular_domain}
\end{figure}
This example models electric field around a cylinder. This case is chosen in order to show the capability of the enriched model to handle curved interfaces.
In this case, the normal direction to the materials interface, and consequently, the orientation of the intra-elemental jump in the gradient vary in the neighboring cut elements.\color{black}

The sketch of the test case is shown in Fig.~\ref{figure:circular_domain} [a]. The electric potentials of $1V$ and $0V$ are considered along the top and bottom edges of the unit square domain, respectively. The entire domain that is enclosed by the circle contains a homogeneous medium with a permittivity 3 times greater than the permittivity of the surrounding part ($\epsilon_1 = 3\epsilon_2$), this medium is located at $x_c=0.25$ and $y_c= 0.75$. The solution field obtained with the proposed E-FEM on an unstructured non-conforming mesh is presented in Fig.~\ref{figure:circular_domain} [c], showing its capability to capture the discontinuity across the curved interface.\\

\begin{figure}
    \centering
                \includegraphics[width=\textwidth]{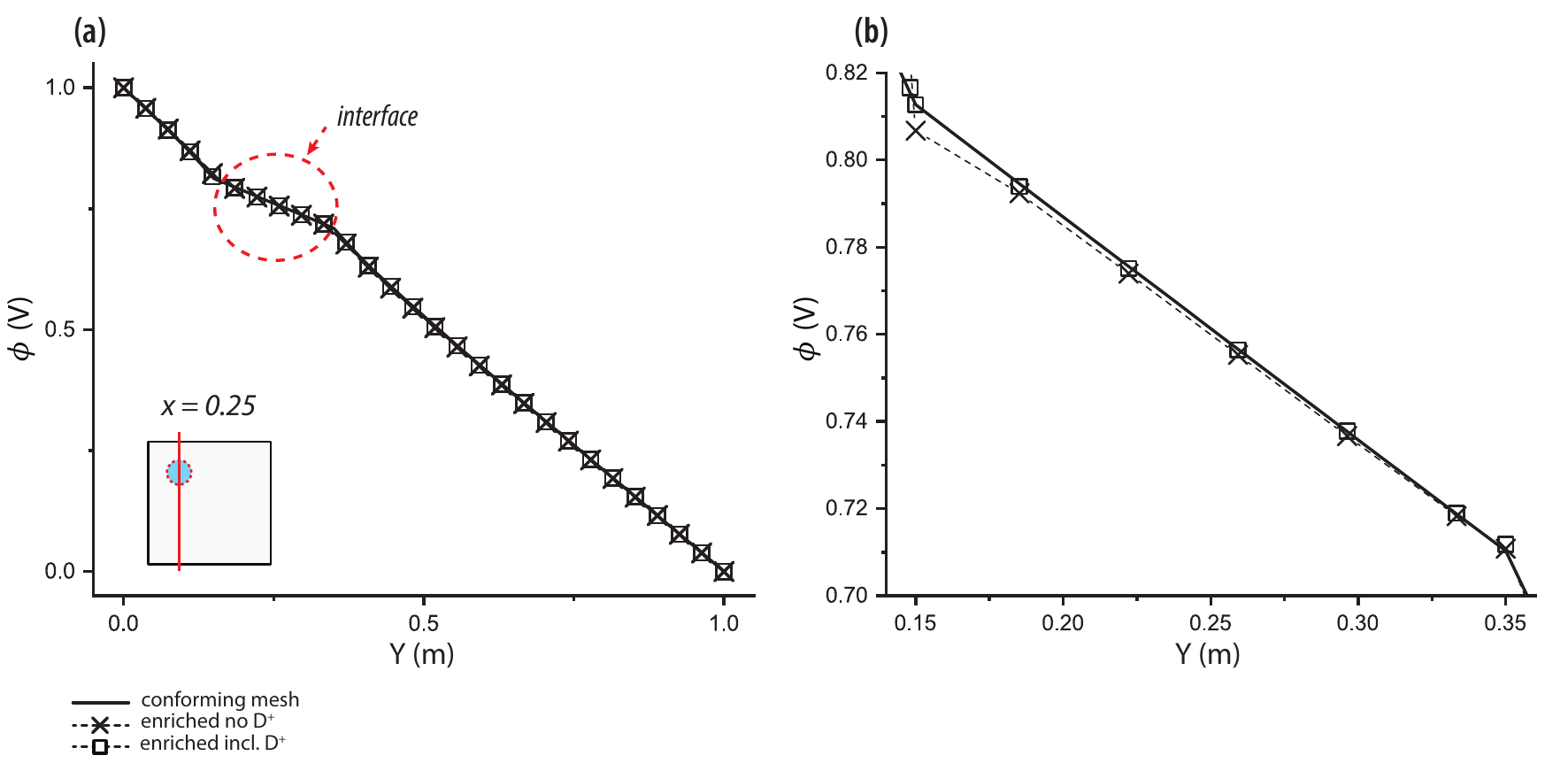}
                \caption{Electric potential distribution around a cylinder. (a) Potential along the vertical cut (x=0.25). (b) Potential at the vicinity of the interface.} 
                \label{figure:circular_potential}
\end{figure}

The solution obtained with standard FEM on a rather fine conforming mesh with 23k elements (h $\approx$ 0.01) is adopted as the reference solution for this test-case. The results for the enriched finite element method are obtained on nonconforming meshes with $h = 0.0375$, $h = 0.01875$, and $h = 0.01$.
In order to be able to precisely visualize the results, after solving the problem, the additional DoFs ($\phi^*$) is further recovered. The results obtained for the coarser mesh ($h = 0.0375$) are compared in Fig.~\ref{figure:circular_potential}, which shows the distribution of the electric potential along a vertical cut ($x=0.25$) that passes though the middle of the cylinder as illustrated in the inset of the figure. By comparison with the reference solution, it is clear that the proposed E-FEM including the D$^+$ term provides quite accurate solutions on the coarse mesh (see Fig.~\ref{figure:circular_potential} [b]).

\begin{figure}
    \centering
                \includegraphics[width=\textwidth]{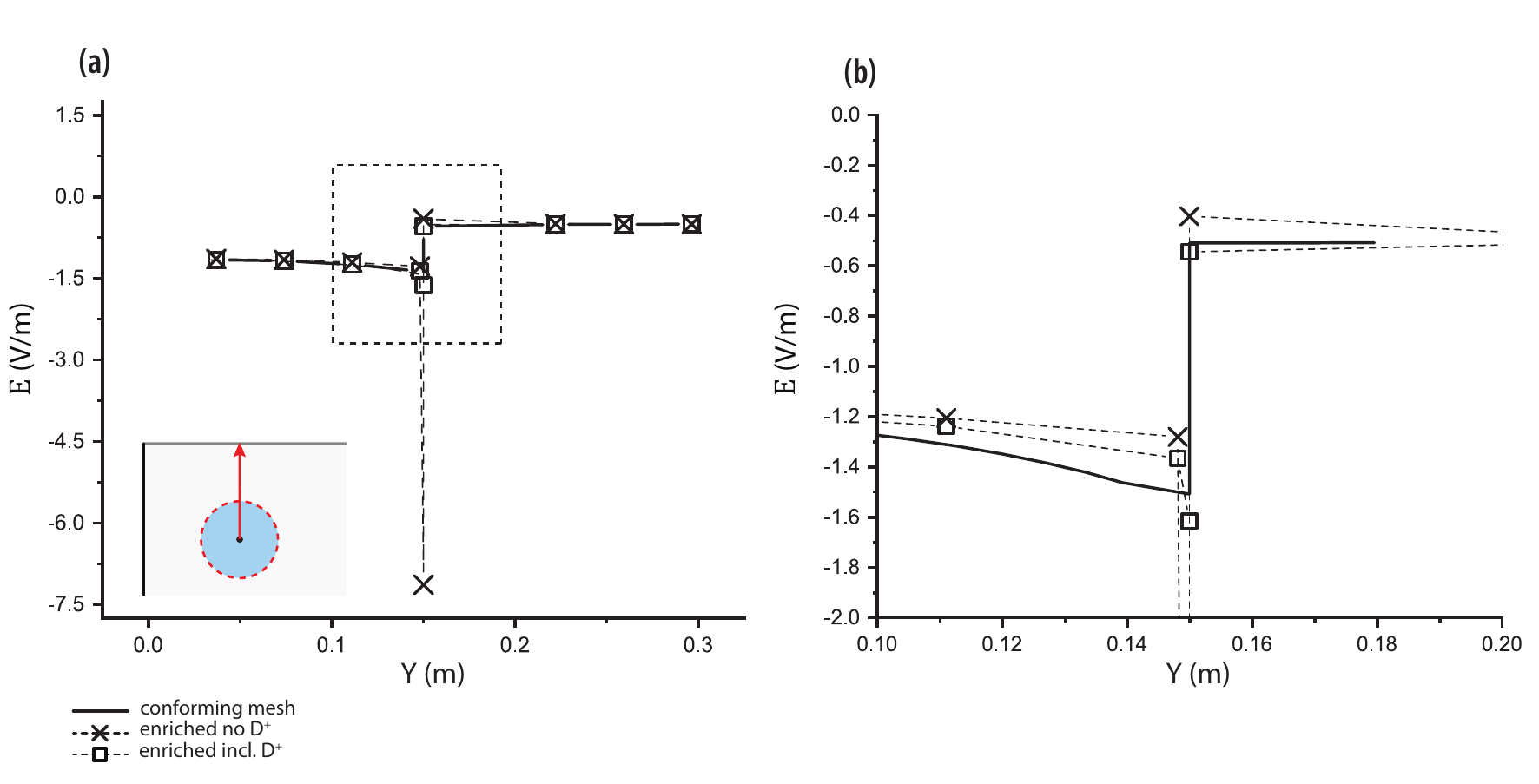}
                \caption{Comparison between reference solution and E-FEM for electric field around cylinder (along a X = 0.25 transect)} 
                \label{figure:circular_E}
\end{figure}

Figure~\ref{figure:circular_E} shows the normal component of the electric field ($E_n$) along the vertical cut for $h = 0.0375$. It is clear that the approximation of the ``jump'' in the electric field at the interface is erroneous when the electric displacement term is neglected. In this case, the error in $E_n$ at the interface is $20\%$ and $38\%$ for medium 1 and 2, respectively. By contrast, using the proposed consistent formulation (E-FEM with the inclusion of D$^+$ term), the errors reduce to $7\%$ in both media. The error reduces by using a more refined mesh as presented in Table~\ref{table:error_circular}. For instance the errors, $E_{n1}$ and $E_{n2}$, decrease to $~3\%$ for the proposed method on a mesh with $h = 0.01875$ (6.3k elements).

\begingroup
\begin{table}
  \centering 
\caption {\label{table:error_circular} Deviations from the reference solution of the computed electric field in the normal direction ($E_n$) using E-FEM with/without the inclusion of D$^+$ term.} 
\begin{tabular}{c cc}
\hline
\multicolumn{1}{c}{element size}   & \multicolumn{2}{c}{error  $|1-E/E^{ref}|$  } \\
\cline{2-3} 
($h$)  & $E_{n1}$ & $E_{n2}$ \\
\midrule
E-FEM no D$^+$ \\
0.0375 & 0.207 & 3.726                    \\
0.01875 & 0.20 & 0.203                    \\
0.01 & 0.32 & 0.026    \\
E-FEM incl. D$^+$\\
0.0375 & 0.069 & 0.072                    \\
0.01875 & 0.034 & 0.036                   \\
0.01 & 0.020 & 0.008    \\ 
\midrule
\end{tabular}
\end{table}
\endgroup

\begin{figure}
    \centering
                \includegraphics[width=\textwidth]{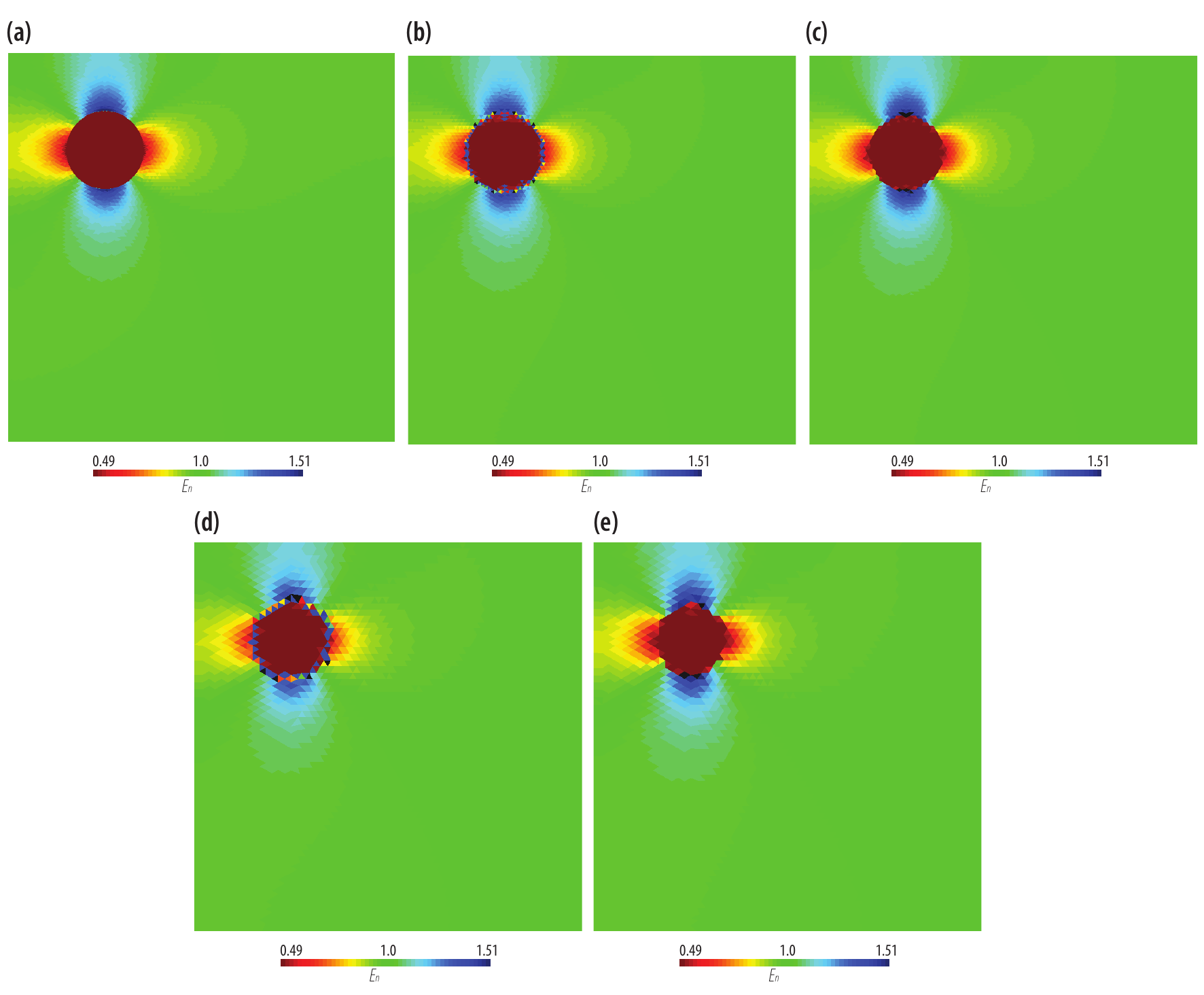}
                \caption{Solution of $E_n$ around cylinder. (a) standard FEM with conforming mesh (h=0.01), (b) E-FEM no D$^+$ (h = 0.01), (c) E-FEM incl. D$^+$ (h = 0.01), (d) E-FEM no D$^+$ (h = 0.01875), and (e) E-FEM incl. D$^+$(h = 0.01875). } 
                \label{figure:numeric_En}
\end{figure}

In order to further analyze the results, the elemental value of the normal component of the electric field is shown in Fig.~\ref{figure:numeric_En} over the whole domain. The spurious values that deteriorate the solution in the vicinity of the interface are resulted from the inconsistency of the formulation in case the inter-elemental electric displacement term is omitted (see Fig.~\ref{figure:numeric_En} [b] and [d]). 
As observed in Fig.~\ref{figure:numeric_En} [c] and [e], the results are significantly improved and the spurious mismatch in the electric field is resolved for both mesh sizes by the inclusion of the inter-elemental D$^+$ term.

In summary, the proposed enriched finite element method provides results for multi-material electrostatic problems with an order of accuracy comparable to that of the mesh-conforming finite element approach.

\subsection{Electric field around a sphere.}

Ultimately we present a 3D simulation considering the electric potential distribution around a sphere. The computational domain is $\Omega = (0,1)\times(0,1)\times(0,1)$. 
The shape is centered at $(0.5,0.5,0,5)$ and the element size is $h\approx0.04$. To apply the electric field, the top and bottom edges located along the y-axis are set to $1V$ and $0V$, respectively.

\begin{figure}
    \centering
                \includegraphics[width=\textwidth]{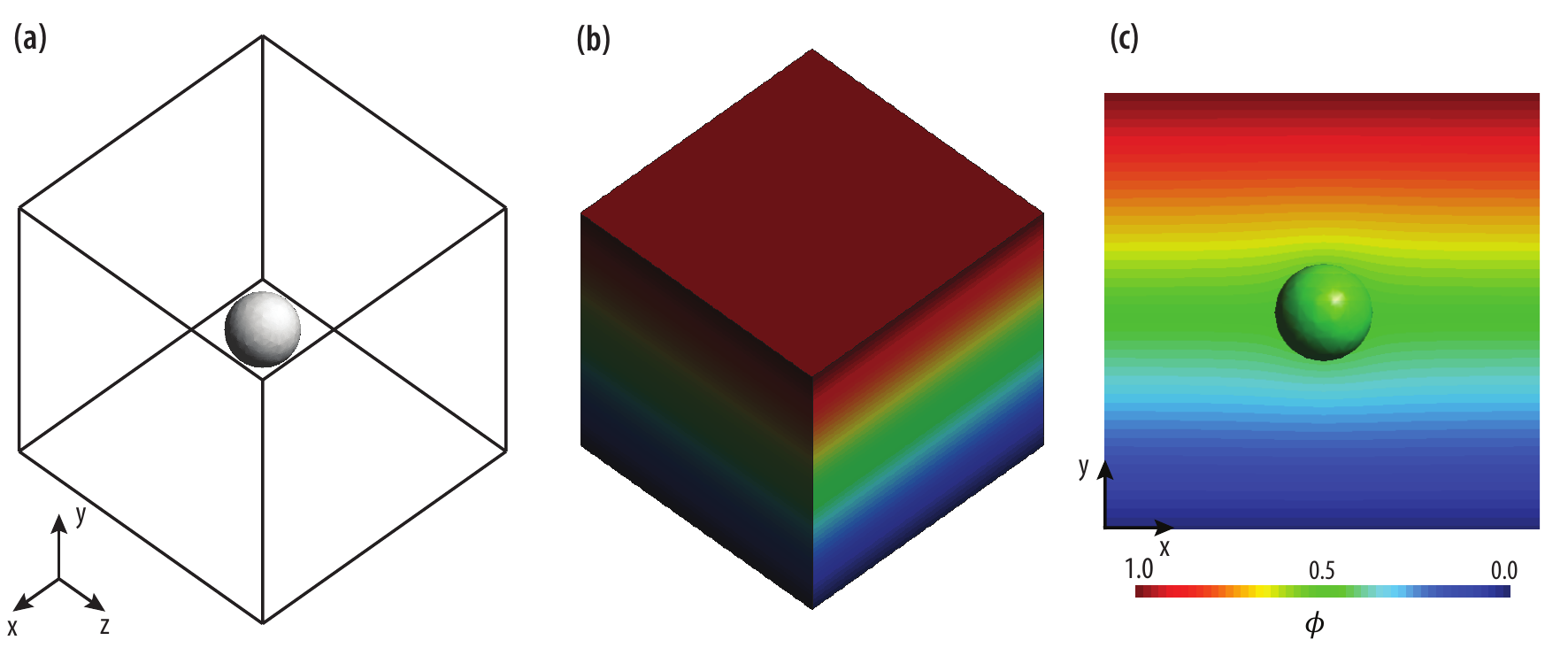}
                \caption{Model problem with a sphere. a) Numerical domain. b) Electric potential distribution along the domain for $Q = 3$. c) Electric potential distribution along the $x-y$ plane.} 
                \label{figure:3d_case}
\end{figure}

A sphere with radius $R_o = 0.1$ is considered (Figure~\ref{figure:3d_case} a)), and the permittivity ratio between the two materials is $Q=3$. For the considered configuration, the analytical solution is given by Eq. \ref{equation:3dpotential} in polar coordinates $(\theta,r)$.
 \begin{equation}
\begin{matrix}

\phi_1 (\theta,r)= \left(\frac{3 r cos(\theta)}{2+Q}\right) \\
\phi_2\theta,r= cos(\theta)\left(r+\frac{1-Q}{2+Q}\frac{R_o^3}{r^2}\right)
\end{matrix}
    \label{equation:3dpotential}
\end{equation}
where $\phi_1$ and $\phi_2$ are the electric potentials in the sphere and the surrounding material, respectively. 

\begin{figure}
    \centering
                \includegraphics{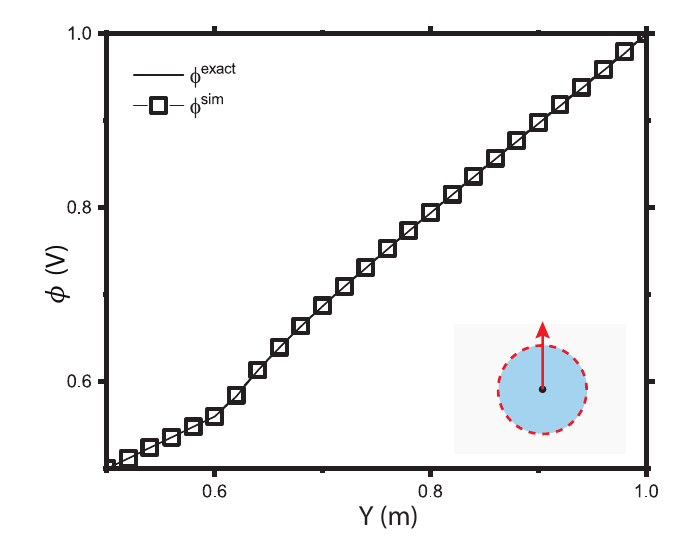}
                \caption{Comparison between reference solution and E-FEM for electric potential around a sphere (along $x-y$ plane)} 
                \label{figure:3d_validation}
\end{figure}

Simulation results are shown in Fig. Fig.~\ref{figure:3d_case} b) and c). For the sake of clear visualization, in the latter sub-figure the electric potential is shown at $x-y$ plane passing though the middle of the sphere. Besides, the numerical results are validated against the exact (analytical) solution in Fig.~\ref{figure:3d_validation}. 
One can see that the numerical solution fits the exact one.

\newpage
\textcolor{Blue}{\section{Summary and conclusions}}
An enriched finite element model for solving electrostatic problems on multi-material domains was developed and implemented in two and three dimensions. The proposed scheme allows using non-conforming mesh at the materials interface. The particular feature of the proposed enriched model is that it uses only one additional degree of freedom per cut element in order to capture the weak discontinuity in the field variable across the material interface. This feature facilitates elemental static condensation of the additional degree of freedom prior to finite element assembly, resulting in an efficient implementation maintaining the graph of the equations system fixed. 
Tangent matrix of the governing equations was formulated, accounting for the element-wise definition of the enriched test functions for obtaining the maximum robustness of the method.

Overall, one can conclude that: 
\begin{itemize}
    \item The proposed E-FEM method accurately represents the discontinuous electric field in multi-material problems, not requiring interface-conforming meshes
    \item  The method is capable of dealing with a wide range of properties ratios (Q), including very large ones ($1/Q\approx0$).
    \item Accounting for inter-elemental electric displacement term in cut elements becomes crucial for small ratios of the material properties.
    \item The proposed scheme is capable of providing accurate results for arbitrary interface orientation, exhibiting nearly second order of accuracy comparable to the results obtained on interface-conforming mesh in all case-studies performed. 
\end{itemize}

Although the methodology is based on multi-material electrostatic problems, it can be used straight-forward in electromagnetic problems. Moreover, the proposed scheme can be easily generalized for any problem represented by the Laplace equation (such as potential flow, heat transfer, diffusion, torsion, among others), characterized by a weak discontinuity of the field variable.   
\newpage
\bibliographystyle{elsarticle-num-names}
\bibliography{sample.bib}
\newpage

\end{document}